\documentclass[a4paper,10pt]{article}
\usepackage{amssymb, amscd, verbatim}
\usepackage[arrow, matrix, curve]{xy}
\usepackage{hyperref}
\usepackage{amsmath}
\usepackage{geometry}

\title{Independence of $\ell$-adic Galois representations over function fields}
\author{Wojciech Gajda \and Sebastian Petersen\footnote{The corresponding author}
}

\begin{document}

\parindent0em
\parskip1em

\newtheorem{leer}{}[section]
\newtheorem{thm}[leer]{Theorem}
\newtheorem{conj}[leer]{Conjecture}
\newtheorem{defi}[leer]{Definition}
\newtheorem{rema}[leer]{Remark} 
\newtheorem{prop}[leer]{Proposition}
\newtheorem{lemm}[leer]{Lemma}
\newtheorem{coro}[leer]{Corollary}
\newtheorem{quest}[leer]{Question}

\newcommand{\ilim}{\mathop{\varinjlim}\limits}
\newcommand{\plim}{\mathop{\varprojlim}\limits}

\newcommand{\OO}{{\cal O}}
\newcommand{\GG}{{\cal G}}
\newcommand{\BB}{{\cal B}}
\newcommand{\KK}{{\cal K}}
\newcommand{\XX}{{\cal X}}
\newcommand{\YY}{{\cal Y}}
\newcommand{\XP}{X_{\ol{P}}}
\newcommand{\Xx}{X_{\ol{\xi}}}

\newcommand{\Gal}{{\mathrm{Gal}}}
\newcommand{\Spec}{\mathrm{Spec}}
\newcommand{\GL}{\mathrm{GL}}
\newcommand{\Hom}{\mathrm{Hom}}
\newcommand{\trdeg}{\mathrm{trdeg}}
\newcommand{\mf}{\mathfrak}
\newcommand{\ol}[1]{\overline{#1}}
\newcommand{\ul}[1]{\underline{#1}}
\newcommand{\Aut}{\mathrm{Aut}}
\newcommand{\End}{\mathrm{End}}
\newcommand{\ke}{\mathrm{ker}}
\newcommand{\im}{\mathrm{im}}
\newcommand{\chara}{\mathrm{char}}
\newcommand{\Hrm}{\mathrm{H}}
\newcommand{\Rrm}{\mathrm{R}}

\def\Pp{{\mathbb{P}}}
\def\Ff{{\mathbb{F}}}
\def\Rr{{\mathbb{R}}}
\def\Cc{{\mathbb{C}}}
\def\Qq{{\mathbb{Q}}}
\def\Zz{{\mathbb{Z}}}
\def\Aa{{\mathbb{A}}}
\def\Nn{{\mathbb{N}}}
\def\Gg{{\mathbb{G}}}
\def\Ll{{\mathbb{L}}}

\maketitle

\begin{abstract}
Let $K$ be a finitely generated extension of $\Qq$. We consider the family of $\ell$-adic representations ($\ell$ varies through the set of all prime numbers) 
of the absolute Galois group of $K$, attached to $\ell$-adic cohomology of a separated scheme of finite type over $K$. We prove that the fields cut out from the algebraic closure of $K$ by the kernels of the representations of the family are linearly disjoint over a finite extension of K. 
This gives a positive answer to a question of Serre. 
\footnotetext{%2000 MSC numbers
\textit{\bf 2010 MSC:}
11G10, 14F20.
}
\footnotetext{%key words and phrases
\textit{\bf Key words:}
Galois representation, \'etale cohomology, abelian variety, finitely generated field.
}
\end{abstract}

%\tableofcontents

\section{Introduction}
\newcommand{\famrho}{(\rho_i)_{i\in I}}
\newcommand{\linL}{{\ell\in\Ll}}
Let $\Gamma$ be a profinite group and $(\Gamma_i)_{i\in I}$ a family of  
groups. For every $i$ let $\rho_i\colon \Gamma\to \Gamma_i$ be a 
homomorphism. Following Serre (cf. \cite[p. 1]{bible}),
we shall say that the family $\famrho$ is
{\em independent}, provided the homomorphism
$$\Gamma\buildrel \rho \over\longrightarrow \prod_{i\in I} \rho_i(\Gamma)$$
induced by the $\rho_i$ is surjective. 
Let $\Gamma'\subset \Gamma$ be a closed subgroup. We 
call the family $\famrho$ {\em independent over $\Gamma'$}, if
$\rho(\Gamma')=\prod_{i\in I}\rho_i(\Gamma')$. Finally we call the
family $\famrho$ {\em almost independent}, if there exists an open subgroup 
$\Gamma'\subset \Gamma$, such that $\famrho$ is independent over $\Gamma'$.
Of particular interest is the special case where $\Gamma=\Gal_K$ is the absolute Galois 
group of a field $K$, and $(\rho_\ell)_{\ell\in\Ll}$ is a family of $\ell$-adic 
representations of $\Gal_K$, indexed by the set $\Ll$ of all prime numbers.
%=G(\Omega/K)$ is 
%the Galois group of a Galois extension $\Omega/K$ of fields, and 
%$\rho_i: G(\Omega/K)\to \Gamma_i$ is a homomorphism into some 
%profinite group $\Gamma_i$ for every $i\in I$. In that case we let $K_i$ be
%the fixed field of the kernel of $\rho_i$. Then $K_i/K$ is a Galois extension
%and $\rho_i$ induces an isomorphims $G(K_i/K)\cong \rho_i(G(\Omega/K))$ for 
%every $i\in I$. Furthermore it is not difficult to see that 
%the family $(\rho_i)_{i\in I}$ is independent if and only if the 
%family $(K_i)_{i\in I}$ of fields is linearly disjoint over $K$.   

Important examples of such families of representations arise as follows:
Let $K$ be a field of characteristic zero and let $X/K$ be a separated $K$-scheme of finite type.
Denote by $\widetilde{K}$ an algebraic closure of $K$. For every $\ell\in\Ll$ and every $q\ge 0$ we consider
the representation of the absolute Galois group $\Gal(\widetilde{K}/K)$
$$\begin{xy}
\xymatrix{
\rho_{\ell, X}^{(q)}\colon \Gal(\widetilde{K}/K)\ar[r] &
\Aut_{\Qq_\ell}(\Hrm^q(X_{\widetilde{K}}, \Qq_\ell))
}
\end{xy}$$
afforded by the \'etale cohomology group $\Hrm^q(X_{\widetilde{K}}, \Qq_\ell)$, and also the
representation
$$\begin{xy}
\xymatrix{
\rho_{\ell, X, c}^{(q)}\colon \Gal(\widetilde{K}/K)\ar[r] &
\Aut_{\Qq_\ell}(\Hrm^q_\mathrm{c}(X_{\widetilde{K}}, \Qq_\ell))
}
\end{xy}$$
afforded by the \'etale cohomology group with compact support $\Hrm^q_\mathrm{c}(X_{\widetilde{K}}, \Qq_\ell)$.
One can wonder in which circumstances the families $(\rho_{\ell, X}^{(q)})_{\ell\in\Ll}$
and $(\rho_{\ell, X, c}^{(q)})_{\ell\in\Ll}$ are almost independent.

In the recent paper \cite{bible} Serre considered the special case where $K$ is a {number field}. He
proved a general independence criterion for certain families of $\ell$-adic representations over a number field 
(cf. \cite[Section 2, Th\'eor\`em 1]{bible}),
and used this criterion together with results of Katz-Laumon and of Berthelot (cf. \cite{illusie})
in order to prove the following 
Theorem (cf. \cite[Section 3]{bible}). 

{\em Let $K$ be a number field and $X/K$ a 
separated scheme of finite type. Then the families of representations $(\rho_{\ell, X}^{(q)})_{\ell\in\Ll}$
and $(\rho_{\ell, X, c}^{(q)})_{\ell\in\Ll}$ are almost independent.} 

The special case of an abelian variety $X$ over a number field $K$ had been dealt with earlier in a letter from Serre to Ribet (cf. \cite{serretoribet2}).
In \cite[p. 4]{bible} Serre asks the following question.

{\em Does this theorem remain true, if one replaces the number field $K$ by a 
finitely generated transcendental extension $K$ of $\Qq$?}  

This kind of problem also shows up in Serre's article \cite[10.1]{serre1994} and in 
Illusie's manuscript \cite{illusie}.
The aim of our paper is to answer this question affirmatively.
In order to do this we prove an independence criterion for families of 
$\ell$-adic representations of the \'etale fundamental group $\pi_1(S)$ of a normal $\Qq$-variety $S$
(cf. Theorem \ref{crit} below). This criterion allows us to reduce the proof of the following 
Theorem \ref{main2} to the 
number field case, where it is known to hold true thanks to the theorem  
of Serre (cf. \cite{bible}) mentioned above. We do take Tate twists into account. 
For every $\ell\in \Ll$ we 
denote by $\varepsilon_\ell\colon \Gal_K\to \Aut_{\Qq_\ell}((\plim_{i\in \Nn} \mu_{\ell^i})\otimes \Qq_\ell)\subset\Qq_\ell^\times$
the cyclotomic character, by $\varepsilon_\ell^{\otimes-1}$ its contragredient and define for every $d\in \Zz$
$$\rho^{(q)}_{\ell, X}(d):=\rho^{(q)}_{\ell, X}\otimes \varepsilon_\ell^{\otimes d}\ \ \mbox{and}\ \ 
\rho^{(q)}_{\ell, X, \mathrm{c}}(d):=\rho^{(q)}_{\ell, X, \mathrm{c}}\otimes \varepsilon_\ell^{\otimes d}.$$

\begin{thm} Let $K$ be a finitely generated extension of $\Qq$. Let $X/K$ be a separated scheme of 
finite type. Then for every $q\in\Nn$ and every $d\in\Zz$  
the families $(\rho^{(q)}_{\ell, X}(d))_\linL$ and
$(\rho^{(q)}_{\ell, X, c}(d))_\linL$ of representations 
of $\Gal_K$ are almost independent. \label{main2}
\end{thm}

Note that outside certain special cases it is not known whether 
the representations occuring in Theorem \ref{main2} are semisimple.
Hence we cannot 
use techniques like the semisimple approximation
of monodromy groups in the proof of Theorem \ref{main2}. 

Theorem \ref{main2} has an important consequence for the arithmetic of abelian varieties.
Let $A/K$ be an abelian variety.
For every $\ell\in\Ll$ consider the Tate module $T_\ell(A):=\plim_{i} A(\widetilde{K})[\ell^i]$, 
define $V_\ell(A):=T_\ell(A)\otimes_{\Zz_\ell} \Qq_\ell$ and let 
$$\begin{xy}
\xymatrix{
\eta_{\ell, A}\colon \Gal(\widetilde{K}/K)\ar[r] & \Aut_{\Qq_\ell}(V_\ell(A))
}
\end{xy}$$
be the $\ell$-adic representation attached to $A$. Then the $\Qq_\ell[\Gal_K]$-modules $V_\ell(A)$ and
$\Hrm^1(A^\vee_{\widetilde{K}}, \Qq_\ell(1))$ are isomorphic, i.e. the representation
$\eta_{\ell, A}$ is isomorphic to $\rho_{\ell, A^\vee}(1)$. 
Hence Theorem \ref{main2} implies that
the family $(\eta_{\ell, A})_{\ell\in\Ll}$ is almost independent. Denote by
$K(A[\ell^\infty])$ the fixed field in $\widetilde{K}$ of the kernel of $\eta_{\ell, A}$. Then
$K(A[\ell^\infty])$ is the field obtained from $K$ by adjoining the coordinates of the
$\ell$-power division points in $A(\widetilde{K})$. Using Remark \ref{indeprem} below we
see that Theorem \ref{main2} has the following Corollary.

\begin{coro} Let $K$ be a finitely generated extension of $\Qq$ and $A/K$ an abelian variety.
Then there is a finite extension $E/K$ such that 
the family $(EK(A[\ell^\infty]))_{\ell\in\Ll}$ is linearly disjoint over $E$.\label{main1}
\end{coro}

This paper carries an appendix with a more elementary proof of this Corollary, which is based on 
our Theorem \ref{crit} below, but avoiding use of \'etale cohomology.

\begin{center} {\bf Notation and Preliminaries} \end{center}

For a field $K$ fix an algebraic closure $\widetilde{K}$
%, let $K_s$ be the separable closure of $K$ in $\widetilde{K}$ 
and denote by $\Gal_K$ the absolute 
Galois group of $K$. %We say that $K$ is {\em finitely generated}, if it is finitely 
%generated over its prime field. 
We denote by $\Ll$ the set of all prime numbers.

Let $S$ be a scheme and $s\in S$ a point (in the underlying topological space). Then
$k(s)$ denotes the residue field at $s$. A {\em geometric point} of $S$ is a 
morphism $\ol{s}\colon \Spec(\Omega)\to S$ where $\Omega$ is an algebraically closed field. To 
give such a geometric point $\ol{s}$ is equivalent to giving a pair $(s, i)$ consisting
of a usual point $s\in S$ and an embedding $i\colon k(s)\to \Omega$. We then let $k(\ol{s})$ be
the algebraic closure of $i(k(s))$ in $\Omega$. Now assume $S$ is an integral scheme and 
let $K$ be its function field. Then we view $S$ as equipped with the geometric
generic point $\Spec(\widetilde{K})\to S$ and denote by $\pi_1(S)$ the \'etale fundamental
group of $S$ with respect to this geometric point.
A {\em variety} $S$ over a field $F$ is an integral separated  
$F$-scheme of finite type. %In that case we often write $F(s)$ for the residue 
%field at $s\in S$ and $F(S)$ for the
%function field of $S$. 

Now let $S$ be a connected normal scheme with function field $K$. Assume for 
simplicity that $\chara(K)=0$. 
If $E/K$ is a an algebraic field extension, then $S^{(E)}$ 
denotes the normalization of $S$ in $E$ (cf. \cite[6.3]{EGAII}). This notation is used throughout this 
manuscript. The canonical morphism $S^{(E)}\to S$ is universally closed and 
surjective. (This follows from the going-up theorem, cf. \cite[6.1.10]{EGAII}.)
If $E/K$ is a finite extension,
then $S^{(E)}\to S$ is a finite morphism (cf. \cite[Proposition I.1.1]{milne}).
We shall say that an algebraic extension $E/K$ is {\em unramified along $S$}, provided 
the morphism $S^{(E')}\to S$ is 
\'etale for every {\em finite} extension $E'/K$ contained in $E$. We denote by 
$K_{S, \mathrm{nr}}$ the maximal extension of $K$ inside $\widetilde{K}$ which is unramified 
along $S$, and by $S_{\mathrm{nr}}$ the normalization of $S$ in $K_{S, \mathrm{nr}}$. One can then 
identify $\pi_1(S)$ with $\Gal(K_{S, \mathrm{nr}}/K)$.  
Let $E/K$ be a Galois extension.
If $P\in S$ is a closed point and $\hat{P}$ is a point in $S^{(E)}$ above $P$, then
we define $D_{E/K}(\hat{P})\subset \Gal(E/K)$ to be the decomposition group
of $\hat{P}$, i.e. the stabilizer of $\hat{P}$ under the action of $\Gal(E/K)$. 
Then $k(\hat{P})/k(P)$ is Galois and the restriction map
$$\begin{xy}
\xymatrix{
r_{E/K, \hat{P}}\colon D_{E/K}(\hat{P})\ar[r] & \Gal(k(\hat{P})/k(P))
}
\end{xy}$$
is an epimorphism. To see this apply \cite[Proposition 1.1, p. 106]{SGA1} for the case $[E:K]<\infty$ and use
a limit argument. If $E/K$ is unramified along $S$, then $r_{E/K, \hat{P}}$ is bijective.

Let $\GG/S$ be a finite \'etale group scheme and choose a closed
point $\hat{P}\in S_{\mathrm{nr}}$. Then there is a finite
extension $E/K$ in $K_{S, \mathrm{nr}}$ such that $\GG\times_S S^{(E)}$ is
a constant group scheme over $S^{(E)}$. In particular, the action of
$\Gal_K$ on $\GG(\widetilde{K})$ factors through $\Gal(K_{S, \mathrm{nr}}/K)$, and the canonical 
evaluation maps 
$$\begin{xy}
\xymatrix{
\GG(S^{(E)})\ar[r] & \GG(E)=\GG(\widetilde{K}) &\mbox{and}& \GG(S^{(E)})\ar[r] & \GG(k(\hat{P}))
}
\end{xy}$$
are bijective. The composite isomorphism 
$$sp_{\GG, \hat{P}}\colon \GG(\widetilde{K})\cong \GG(k(\hat{P}))$$
is called the {\em cospecialization map}. This map is equivariant in the sense that 
$sp_{\GG, \hat{P}}(\sigma (x))=r_{E/K, P}(\sigma)(sp_{\GG, \hat{P}}(x))$ for all 
$x\in\GG(E)$ and all $\sigma\in D_{E/K}(\hat{P})$.   
\label{preliminary}

\section{Finiteness properties of Jordan extensions}
Let $E/K$ be an algebraic field extension and $d\in \Nn$. We will say that
$E/K$ is {\em $d$-flat}, if $E$ is a compositum of (finitely or infinitely many) 
Galois extensions of $K$, each of degree $\le d$. In particular every $d$-flat extension is Galois.
We call the extension $E/K$
{\em $d$-Jordanian}, if $E/K$ is a (possibly infinite) abelian extension of a $d$-flat extension.
The $1$-Jordanian extensions of $K$ are hence just the abelian extensions of $K$.
If $K$ is a number field and $E/K$ is a $d$-Jordanian extension of $K$ which is
everywhere unramified, then $E/K$ is finite. This has been shown by Serre in \cite[Th\'eor\`eme 2]{bible},
making use of the Hermite-Minkowski theorem and the finiteness of the Hilbert class field.
The aim of this section is to derive a similar finiteness property for $d$-Jordanian
extensions of function fields over $\Qq$. In Lemmata \ref{mwlemm}, \ref{curves} and \ref{katzlanglemm}
we follow closely the paper \cite{katzlang} of Katz and Lang on geometric class field theory, giving
complete details for the convenience of the reader.

If $E$ is any extension field of $\Qq$, then we denote by $\kappa_E$ the 
algebraic closure of $\Qq$ in $E$, 
$$\kappa_E:=\{x\in E\,: x\ \mbox{is algebraic over}\ \Qq\},$$
and we call $\kappa_E$ the {\em constant field} of $E$. 

\begin{rema} Let $K$ be a finitely generated extension of $\Qq$.
Let $E/K$ be an algebraic extension. Then there is a diagram of fields:
$$\begin{xy}
  \xymatrix{
  & E\ar@{-}[r]\ar@{-}[d] & \widetilde{\Qq}E \ar@{-}[d]\\
\kappa_E\ar@{-}[r]\ar@{-}[d]  & \kappa_EK \ar@{-}[r] \ar@{-}[d] & \widetilde{\Qq}K \\
\kappa_K\ar@{-}[r] & K & 
}
\end{xy}$$
The field $\kappa_K$ is a number field and $\kappa_E/\kappa_K$
is an algebraic extension. We say that $E/K$ is a {\em constant field extension}, if
$\kappa_E K=E$. If $E/K$ is Galois, then $\kappa_E/\kappa_K$,
$\kappa_EK/K$ and $\widetilde{\Qq}E/\widetilde{\Qq}K$ are Galois as well, and
the restriction maps $\Gal(\widetilde{\Qq}E/\widetilde{\Qq}K)\to \Gal(E/\kappa_E K)$ and
$\Gal(\kappa_E K/K)\to \Gal(\kappa_E/\kappa_K)$ are both bijective.\label{elrema}
\end{rema}

The aim of this section is to prove the following Proposition.

\begin{prop} Let $S/\Qq$ be a normal variety with function field $K$. Let $d\in \Nn$. 
Let $E/K$ be a $d$-Jordanian 
extension which is unramified along $S$. Then $E/\kappa_E K$ is a finite extension. \label{jordan}
\end{prop}

Note that in the situation of Proposition \ref{jordan} the extension $\kappa_E/\kappa_K$ 
may well be infinite algebraic.
The proof occupies the rest of this section.

\begin{lemm} Let $S/\Qq$ be a normal variety with function field $K$. 
Let $d\in \Nn$. Let $E/K$ be a $d$-flat 
extension which is unramified along $S$. Then $E/\kappa_E K$ is finite and
$\Gal(\kappa_E/\kappa_K)$ is a (possibly infinite) group of exponent $\le d$. \label{fin1}
\end{lemm}

{\em Proof.} There is a sequence $(K_i)_{i\in I}$ of intermediate fields of $E/K$ 
such that each $K_i/K$ is Galois with $[K_i:K]\le d$ and $E=\prod_{i\in I} K_i$. 
Hence $\Gal(E/K)$ is a closed subgroup of $\prod_{i\in I} \Gal(K_i/K)$. By Remark \ref{elrema}
$\Gal(\kappa_E/\kappa_K)$ is a quotient of $\Gal(E/K)$, hence  
$\Gal(\kappa_E/\kappa_K)$ has exponent $\le d$.
Again by Remark \ref{elrema} it is now enough to show that $\widetilde{\Qq}E/\widetilde{\Qq}K$
is finite. The Galois group $\Gal(\widetilde{\Qq}E/\widetilde{\Qq}K)$ 
is a quotient of $\pi_1(S_{\widetilde{\Qq}})$, and
$\pi_1(S_{\widetilde{\Qq}})$ is topologically finitely generated (cf. \cite[II.2.3.1]{SGA7}). Hence there 
are only finitely many intermediate fields $L$ of $\widetilde{\Qq}E/\widetilde{\Qq}K$ with
$[L:\widetilde{\Qq}K]\le d$ (cf. \cite[16.10.2]{friedjarden}). This implies that $\widetilde{\Qq}E/\widetilde{\Qq}K$ 
is finite.\hfill $\Box$

\begin{lemm} \label{exprime} Let $K$ be a finitely generated extension of $\Qq$. 
Let $E/K$ be a (possibly infinite) Galois extension. Assume that $\Gal(E/K)$ has
finite exponent. Let $X=(X_1,\cdots, X_n)$ be a transcendence base 
of $K/\Qq$ and $R$ the integral closure of $\Zz[X]$ in $E$. 
\begin{enumerate}
\item[a)] The residue field $k(\mf m)=R/\mf m$ is {\em finite} for
every maximal ideal $\mf m$ of $R$.
\item[b)] For every non-zero element $f\in R$ there exist two 
maximal ideals $\mf m_1$ and $\mf m_2$ of $R$ such that 
$f\notin \mf m_1$ and $f\notin \mf m_2$ and $\chara(k(\mf m_1))\neq
\chara(k(\mf m_2))$. 
\end{enumerate} 
\end{lemm} 

{\em Proof.} Let $R'$ be the integral closure of $\Zz[X]$ in $K$. 
Let $\mf m$ be a maximal ideal of $R$. Define $\mf m':=\mf m\cap R'$ and 
$\mf p=\mf m\cap \Zz[X]$. There are diagrams of fields and residue fields 
$$\begin{xy}
\xymatrix{
\Qq(X)\ar@{-}[r] & K \ar@{-}[r] & E & \mbox{and}& k(\mf p)\ar@{-}[r] & k(\mf m') \ar@{-}[r] & k(\mf m).
}
\end{xy}$$
By the going-up theorem $\mf p$ is a
maximal ideal of $\Zz[X]$, and $k(\mf p)=\Zz[X]/\mf p$ is a {\em finite}
field. Furthermore $R'$ is a finitely generated $\Zz[X]$-module
(cf. \cite[Prop. I.1.1]{milne}). This implies that 
$k(\mf m')$ is a finite field. The extension $k(\mf m)/k(\mf m')$ is
Galois and the Galois group $G:=\Gal(k(\mf m)/k(\mf m'))$ is a subquotient of 
$\Gal(E/K)$. Hence $G$ is of finite exponent. On the other hand $G$ must be
procyclic, because it is a quotient of the Galois group $\hat{\Zz}$ of the
finite field $k(\mf m')$. It follows that $G$ is finite and that 
$k(\mf m)$ is a finite field. This finishes the proof of part a).

Now let $f\in R$ be a nonzero element. The canonical
morphism $p: \Spec(R)\to\Spec(\Zz[X])$ is closed (cf. \cite[6.1.10]{EGAII}), 
hence $p(V(f))$ is a closed subset of $\Spec(\Zz[X])$. It is also a proper subset 
of $\Spec(\Zz[X])$. 
It follows that there is a non-zero polynomial $g\in \Zz[X]$ such that 
$D(g)\cap p(V(f))=\emptyset$. Choose $a\in
\Zz^n$ with $g(a)\neq 0$. Then choose distinct prime numbers $p_1\neq p_2$ 
not dividing $g(a)$. For $i\in\{1, 2\}$ consider 
the maximal ideal $\mf p_i=(p_i, X-a_1,\cdots, X-a_n)$ of $\Zz[X]$.
Then $\mf p_1, \mf p_2\in D(g)$. Finally let 
$\mf m_1$ and $\mf m_2$ be prime ideals of $R$ such that $p(\mf m_i)=\mf p_i$ for
$i\in\{1, 2\}$. Then $\mf m_1$ and $\mf m_2$ have the desired 
properties.\hfill $\Box$

We now show that a weak form of the Mordell-Weil theorem holds true over finitely generated
extensions of fields like the field $\kappa_E$ occuring in Lemma \ref{fin1}. If $B$ is a semiabelian variety
over a field $K$, then we define $T(B)=\prod_{\ell\in \Ll} T_\ell(B)$
and $T(B)_{\neq p}:=\prod_{\ell\in \Ll\setminus\{p\}} T_\ell(B)$ (for $p\in\Ll$), 
where $T_\ell(B)=\plim_{i\in\Nn} B(\widetilde{K})[\ell^i]$ is the Tate module of $B$ for every $\ell\in \Ll$.
If $M$ is a compact topological $\Gal_K$-module, then we define the {\em module of coinvariants}
$M_{\Gal_K}$ of $M$ to be the largest Hausdorff quotient of $M$ on which $\Gal_K$ acts trivially.
% that is:  $$M_{\Gal_K}=M/\overline{\langle \sigma x-x|\sigma\in\Gal_K, x\in M\rangle}.$$ 
%Thus $M_{\Gal_K}$ is the maximal Hausdorff quotient of
%$M$ on which $\Gal_K$ acts trivially.
 
\begin{lemm} Let $K$ be a finitely generated extension of $\Qq$. Let 
$E/K$ be a Galois extension. Assume that $\Gal(E/K)$ has
finite exponent. Let $B/E$ be a semiabelian variety. Then
$T(B)_{\Gal_E}$ is finite. \label{mwlemm}
\end{lemm} 

{\em Proof.} Let $E'/E$ be a finite extension over which the torus part of $B$ 
splits. Then there exists a finite
Galois extension $L/K$ such that $LE\supset E'$, and $\Gal(LE/K)$ has 
finite exponent again. The group $T(B)_{\Gal_{E}}$ is a quotient of $T(B)_{\Gal_{LE}}$. 
Hence we may assume right from
the beginning that $B$ is an extension of an abelian variety $A$ by
a split torus $\Gg_{m, E}^d$. 
Then there is an exact sequence of $\Gal_E$-modules
$$\begin{xy}
\xymatrix{
0\ar[r] & T(\Gg_m)^d\ar[r] & T(B)\ar[r] & T(A)\ar[r] & 0
}
\end{xy}$$
As the functor $-_{\Gal_E}$
is right exact, it is enough to prove that $T(A)_{\Gal_E}$ and
$T(\Gg_m)_{\Gal_E}$ are both finite. We may thus assume
that either $B$ is an abelian variety over $E$ (case 1) or 
$B=\Gg_{m, E}$ (case 2). We shall prove the finiteness of 
$T(B)_{\Gal_E}$ in both cases.    

Choose a transcendence base $X=(X_1,\cdots, X_n)$ of $K/\Qq$ and let
$R$ be the integral closure of $\Zz[X]$ in
$E$. In case 1 there is a nonempty open subscheme $U\subset \Spec(R)$ 
such that $B$ extends to an abelian scheme $\BB$ over $U$. In case 2 we
define $U=\Spec(R)$ and put $\BB:=\Gg_{m, U}$. Let $\mf m$ be a  
maximal ideal of $R$ contained in $U$, define $p=\chara(R/\mf m)$, and
denote by $\ol{B}=\BB\times_U \Spec(k(\mf m))$ the special fibre at $\mf m$. 
Let $n$ be a positive integer which is coprime to $p$. Then the
restriction of $\BB[n]$ to $S:=U[1/n]$ is a finite \'etale 
group scheme over $S$ and $\mf m\in S$. Let ${\mf m}_\mathrm{nr}$ be a closed point of $S_{\mathrm{nr}}$ over $\mf m$.
Taking a projective limit over the cospecialization maps (cf. Section \ref{preliminary})
$B[n](\widetilde{E})\cong \ol{B}[n](k({\mf m}_\mathrm{nr}))$, we obtain an isomorphism
$$T(B)_{\neq p}\cong T(\ol{B})_{\neq p},$$
which induces a surjection
$T(\ol{B})_{\neq p, \Gal_{\Ff}}\to T({B})_{\neq p, \Gal_E}$, where 
we have put $\Ff=k(\mf m)$. The field $\Ff$ is {\em finite} by Lemma
\ref{exprime} and $\ol{B}$ is either an abelian variety over $\Ff$ (case 1)
or the multiplicative group scheme over $\Ff$ (case 2). In both cases it is 
known that $T(\ol{B})_{\neq p, \Gal_{\Ff}}$ is finite (cf. \cite[Theorem 1 (ter), p. 299]{katzlang}). 
This shows that 
$T(B)_{\neq p, \Gal_E}$ is finite, whenever there exists a maximal ideal $\mf m$ of
$R$ contained in $U$ with $\chara(k(\mf m))=p$. Now it follows by part b) of Lemma \ref{exprime} that 
there are two different prime numbers $p_1\neq p_2$ such that 
$T(B)_{\neq p_1, \Gal_E}$ and $T(B)_{\neq p_2, \Gal_E}$ are finite, and the 
assertion follows from that.\hfill $\Box$

Let $K_0$ be a field of characteristic zero and $S/K_0$ a normal geometrically irreducible
variety with function field $K$. There is a canonical epimorphism
$p\colon \pi_1(S)\to \Gal_{K_0}$ (with kernel $\pi_1(S_{\widetilde{K_0}})$) and, following 
Katz-Lang (\cite[p. 285]{katzlang}), we define $\KK(S/K_0)$ to be the 
kernel of the map $\pi_1(S)_{\mathrm{ab}}\to \Gal_{K_0, \mathrm{ab}}$ induced by $p$ on the abelianizations.
If we denote by $K_{S, \mathrm{nr}, \mathrm{ab}}$ the maximal abelian extension of $K$ which is 
unramifield along $S$, then there is a diagram of fields  
$$\begin{xy}
  \xymatrix{
  & K_{S, \mathrm{nr}, \mathrm{ab}}\ar@{-}[r]\ar@{-}[d] & \widetilde{K_0}K_{S, \mathrm{nr}, \mathrm{ab}} \ar@{-}[d]\\
K_{0, \mathrm{ab}}\ar@{-}[r]\ar@{-}[d]  & K_{0, \mathrm{ab}}K \ar@{-}[r] \ar@{-}[d] & \widetilde{K_0}K \\
K_0\ar@{-}[r] & K & 
}
\end{xy}$$
(cf. \cite[p. 286]{katzlang}) and the groups $\Gal(K_{S, \mathrm{nr}, \mathrm{ab}}/K_{0, \mathrm{ab}}K)$ 
and $\Gal(\widetilde{K_0}K_{S, \mathrm{nr}, \mathrm{ab}}/\widetilde{K_0}K)$ are both isomorphic to $\KK(S/K_0)$. 
The main result in the paper \cite{katzlang} of Katz and Lang is: If $K_0$ is finitely
generated and $S/K_0$ a smooth geometrically irreducible variety, then $\KK(S/K_0)$ 
is finite. On the other hand, if $K_0$ is algebraically closed and $S/K_0$ is a
smooth proper geometrically irreducible curve of genus $g$, 
then $\KK(S/K_0)\cong \hat{\Zz}^{2g}$ is infinite, unless $g=0$. In order to finish
up the proof of Proposition \ref{jordan} we have to prove the 
finiteness of $\KK(S/K_0)$ in the case of certain algebraic extensions $K_0/\Qq$ (like
the field $\kappa_E$ in Lemma \ref{fin1})
which are {\em not} finitely generated but much smaller than $\widetilde{\Qq}$. 
%The proof is only a slight 
%adaptation of the arguments in \cite{katzlang}, based on the observation in
%Lemma \ref{exprime}. Nevertheless we sketch the argument, because the 
%result is of crucial importance in the next section.

\begin{lemm} Let $K$ be a finitely generated extension of $\Qq$.
Let $E/K$ be a (possibly infinite) Galois extension. Assume that $\Gal(E/K)$ has
finite exponent. Let $C/E$ be a smooth proper geometrically irreducible
curve and $S$ the complement of a divisor $D$ in $C$. Then 
$\KK(S/E)$ is finite.\label{curves}
\end{lemm}

{\em Proof.} There is a finite extension $E'/E$ such that $S$ has an
$E'$-rational point and $D$ is $E'$-rational. There is a finite extension
$E''/E'$ which is Galois over $K$. Then $\Gal(E''/K)$ must have finite
exponent (because $\Gal(E/K)$ and $\Gal(E''/E)$ do). Furthermore
$\KK(S_{E''}/E'')$ surjects onto $\KK(S/E)$ (cf. \cite[Lemma 1, p. 291]{katzlang}).
Hence we may assume from
the beginning that $S$ has an $E$-rational point and $D$ is $E$-rational. 
The generalized Jacobian $J$ of $C$ with respect to the modulus
$D$ is a semiabelian variety. 
(If $S=C$, then $J$ is just the usual Jacobian variety of $C$.) 
Furthermore there is an isomorphism
$$\pi_1(S_{\widetilde{E}})_{\mathrm{ab}}\cong T(J).$$
On the other hand $\pi_1(S_{\widetilde{E}})_{\mathrm{ab}, \Gal_E}$ is isomorphic to $\KK(S/E)$ (cf. \cite[Lemma 1, p. 291]{katzlang}).
Hence it is enough to prove that $T(J)_{\Gal_E}$ is finite. But this 
has already been done in Lemma \ref{mwlemm}.\hfill $\Box$

\begin{lemm} Let $K$ be a finitely generated extension of $\Qq$. Let 
$E/K$ be a (possibly infinite) Galois extension. Assume that $\Gal(E/K)$ has
finite exponent. Let $S/E$ be a normal geometrically irreducible variety.
Then $\KK(S/E)$ is finite.\label{katzlanglemm}
\end{lemm} 

{\em Proof.} There is a finite extension $L/E$ and a sequence of elementary fibrations
in the sense of M. Artin (cf. \cite[Expos\'e XI, 3.1-3.3]{SGA4})
$$\begin{xy}
\xymatrix{
\Spec(L)=U_0\ar@{<-}[r]^{f_1} & U_1 \ar@{<-}[r]^{f_2} & U_2\ar@{<-}[r]^{f_3} &\cdots\ar@{<-}[r]^{f_n} & U_n\subset S_L 
}
\end{xy}
$$
where $U_n$ is a non-empty open subscheme of $S_L$. Then $\dim(U_i)=i$ for all $i$. We may assume that $L/K$ is 
Galois and then $\Gal(L/K)$ is of finite exponent. Let $L_i$ be the 
function field of $U_i$. Then the generic fibre $S_{i+1}:=U_{i+1}\times_{U_i} \Spec(L_i)$
of $f_{i+1}$ is a curve over $L_i$ which is the complement of a divisor in a smooth proper
geometrically irreducible curve $C_{i+1}/L_i$. The extension $L_i/L$ is finitely generated
(of transcendence degree $i$). Hence $L_i=L(u_1,\cdots, u_s)$ for certain elements
$u_1,\cdots, u_s\in L_i$. Let us define $K_i:=K(u_1,\cdots, u_s)$. Then there
is a diagram of fields
$$\begin{xy}
\xymatrix{
K_i\ar@{-}[rr]\ar@{-}[d] & & L_i\ar@{-}[d]\\
K\ar@{-}[r]\ar@{-}[d] & E\ar@{-}[r] & L\\
\Qq &&
}
\end{xy}
$$
such that the vertical extensions are all finitely generated and $L_i=K_iL$.
The extension $L_i/K_i$ is Galois because $L/K$ is Galois, and the restriction
map $\Gal(L_i/K_i)\to\Gal(L/K)$ is injective. Hence $\Gal(L_i/K_i)$ is a group of
finite exponent and $K_i$ is finitely generated. 

Lemma \ref{curves} implies that 
$\KK(S_{i+1}/L_i)$ is finite for every $i\in\{0,\cdots, n-1\}$. By 
\cite[Lemma 2]{katzlang} and \cite[(1.4)]{katzlang} it follows that $\KK(U_n/L)$ is finite. 
Then \cite[Lemma 3]{katzlang} implies that $\KK(S_L/L)$ is finite, and \cite[Lemma 1]{katzlang} 
shows that $\KK(S/E)$ is finite, as desired.\hfill $\Box$. 

{\em Proof of Proposition \ref{jordan}.} Let $S/\Qq$ be a normal variety
with the function field $K$. Let $E/K$ be a $d$-Jordan extension contained in 
the extension $K_{S, \mathrm{nr}}/K$.% of $K$ which is unramified along $S$.
There is a $d$-flat extension $L/K$ in $E$ such that $E/L$ is abelian. 
By Lemma \ref{fin1} $L/\kappa_LK$ is a finite extension. We have the following diagram of fields.
$$\begin{xy}
  \xymatrix{
\kappa_E\ar@{-}[r]\ar@{-}[d] & \kappa_E K\ar@{-}[r]\ar@{-}[d] & \kappa_E L\ar@{-}[r]\ar@{-}[d] & E\ar@{-}[r] & K_{S, \mathrm{nr}} \\
\kappa_L\ar@{-}[r]\ar@{-}[d] & \kappa_L K\ar@{-}[r]\ar@{-}[d] & L & & \\
\kappa_K\ar@{-}[r] & K  
}
\end{xy}$$
Now $S^{(L)}$ is the normalization of the geometrically irreducible 
$\kappa_L$-variety $$S^{(\kappa_LK)}=S\times_{\kappa_K}\Spec(\kappa_L)$$ in the 
{\em finite} extension $L/\kappa_LK$. Hence $S^{(L)}$ is a geometrically 
irreducible variety over $\kappa_L$. (The crucial point is that $S^{(L)}$ is
of finite type over $\kappa_L$.) The extension $E/L$ is abelian and unramified along $S^{(L)}$.
Hence $\Gal(E/\kappa_EL)$ is a quotient of $\KK(S^{(L)}/\kappa_L)$. 
The field $\kappa_K$ is a {\em number field} and $\Gal(\kappa_L/\kappa_K)$ is a
group of exponent $d<\infty$, because it is a quotient of $\Gal(L/K)$ (cf. Remark \ref{elrema}). Hence
Lemma \ref{katzlanglemm} implies that $\KK(S^{(L)}/\kappa_L)$ is finite. 
It follows that $E/\kappa_EL$ is a finite extension. 
Now $\kappa_EL/\kappa_EK$ is finite, because $L/\kappa_L K$ is finite. It follows
that $E/\kappa_E K$ is finite, as desired.\hfill $\Box$

\section{Representations of the fundamental group}
We start this section with two 
remarks and a lemma about families of representations of certain profinite groups.
Then we prove an independence criterion for families
of representations of the \'etale fundamental group $\pi_1(S)$ of a normal $\Qq$-variety
$S$ (cf. Theorem \ref{crit}). This criterion is the technical heart of the paper.  

\begin{rema} Let $K$ be a field, $\Omega/K$ a Galois extension and $I\subset \Nn$. Let $(\Gamma_i)_{i\in I}$ be a family of profinite
groups. For every $i\in I$ let $\rho_i: \Gal(\Omega/K)\to \Gamma_i$ be a continuous homomorphism.
Let $K_i$ be the fixed field of $\ker(\rho_i)$ in $\Omega$.
Then the following conditions are equivalent.
\begin{enumerate}
\item[(i)] The family $(\rho_i)_{i\in I}$ is independent.
\item[(ii)] The family $(K_i)_{i\in I}$ of 
fields is linearly disjoint over $K$.
\item[(iii)] If $s\ge 1$ and $i_1<i_2<\cdots < i_{s+1}$ are elements of $I$, then 
$$K_{i_1}\cdots K_{i_s}\cap K_{i_{s+1}}=K.$$
\label{indeprem}
\end{enumerate}
\end{rema}

{\em Proof.} As the homomorphisms $\rho_i$ induce
isomorphisms $\Gal(K_i/K)\cong \im(\rho_i)$, (i) is satisfied 
if and only if the natural map
$\Gal(\Omega/K)\to \prod_{i\in I} \Gal(K_i/K)$ is surjective, and 
this is in turn equivalent to (ii) (cf. \cite[2.5.6]{friedjarden}).
It is well-known that (ii) is equivalent 
to (iii) (cf. \cite[p. 36]{friedjarden}).\hfill $\Box$

\begin{rema} Let $\Gamma$ be a profinite group and $n\in \Nn$. 
For every $\ell\in\Ll$ let $\Gamma_\ell$ be a profinite group and $\rho_\ell\colon \Gamma\to \Gamma_\ell$
a continuous homomorphism. Assume that for every $\ell\in \Ll$ there is an integer $n\in \Nn$ such that $\Gamma_\ell$
is isomorphic to a subquotient of $\GL_n(\Zz_\ell)$.  
\begin{enumerate}
\item[a)] Let $\Gamma'\subset \Gamma$ be an open subgroup.
If the family $(\rho_\ell)_{\ell\in \Ll}$ is independent, 
then there is a finite subset $I\subset \Ll$ such that 
the family $(\rho_\ell)_{\ell\in \Ll\smallsetminus I}$ is independent
over $\Gamma'$.
\item[b)] The following conditions (i) and (ii) are equivalent.
\begin{enumerate}
\item[(i)] The family $(\rho_\ell)_{\ell\in \Ll}$ is almost 
independent.
\item[(ii)] There exists a finite subset $I\subset \Ll$ such
that $(\rho_\ell)_{\ell\in \Ll\smallsetminus I}$ is almost independent.\label{indeprema}
\end{enumerate}
\end{enumerate}
\end{rema}

{\em Proof.} Let $\rho\colon \Gamma\to \prod_{\ell\in\Ll} \Gamma_\ell$ be
the homomorphism induced by the $\rho_\ell$. 
To prove a) assume that $\rho(\Gamma)=\prod_{\ell\in\Ll} \rho_\ell(\Gamma)$. The subgroup
$\rho(\Gamma')$ is open in $\prod_{\ell\in\Ll} \rho_\ell(\Gamma)$, because a
surjective homomorphism of profinite groups is open (cf. \cite[p. 5]{friedjarden}). 
It follows from the definition of the product topology that 
there is a finite subset $I\subset \Ll$ such that 
$\rho(\Gamma')\supset \prod_{\ell\in I}\{1\}\times \prod_{\ell\in\Ll\smallsetminus I} \rho_\ell(\Gamma)$. This
implies that $(\rho_\ell)_{\ell\in\Ll\smallsetminus I}$ is independent over $\Gamma'$ and finishes the proof of part a). 
For part b) see \cite[Lemme 3]{bible}.\hfill $\Box$   

Let $K$ be a field, $n\in\Nn$ and $\Omega/K$ a fixed Galois extension.
For every $\ell\in\Ll$ let $\Gamma_\ell$ be a profinite group and $\rho_\ell\colon \Gal(\Omega/K)\to \Gamma_\ell$
a continuous homomorphism. Assume that $\Gamma_\ell$ is isomorphic to a subquotient of 
$\GL_n(\Zz_\ell)$ for every $\ell\in \Ll$. Denote by $K_\ell$ the fixed field
in $\Omega$ of the kernel of $\rho_\ell$. Then $K_\ell$ is a Galois 
extension of $K$ and $\rho_\ell$ induces an isomorphism
$\Gal(K_\ell/K)\cong \rho_\ell(\Gal(\Omega/K))$.
For every extension $E/K$ contained in $\Omega$ and
every $\ell\in\Ll$ we
define $G_{\ell, E}:=\rho_\ell(\Gal(\Omega/E))$ and $E_\ell:=EK_\ell$. Then
$G_{\ell, E}$ is isomorphic to a subquotient of $\GL_n(\Zz_\ell)$ and
$\rho_\ell$ induces an isomorphism
$$\Gal(E_\ell/E)\cong G_{\ell, E}.$$ Furthermore we define
$G_{\ell, E}^+$ to be the subgroup of $G_{\ell, E}$
generated by its $\ell$-Sylow subgroups. Then $G_{\ell, E}^+$
is normal in $G_{\ell, E}$. Finally we
let $E_\ell^+$ be the fixed field of $\rho_\ell^{-1}(G_{\ell, E}^+)\cap\Gal(\Omega/E)$.
Then $E_\ell^+$ is an intermediate field of $E_\ell/E$ which is
Galois over $E$, the group $\Gal(E_\ell/E_\ell^+)$ is isomorphic to $G_{\ell, E}^+$ 
and $\Gal(E_\ell^+/E)$ is isomorphic to
$G_{\ell, E}/G_{\ell, E}^+$. 
 
\begin{lemm} Let $E/K$ be a Galois extension contained in $\widetilde{K}$ and
let $\ell\in\Ll$.\label{lemm1}
\begin{enumerate}
\item[a)] The extension $E_\ell^+/E$ is a {\em finite} Galois extension, and 
$\Gal(E_\ell^+/E)$ is isomorphic to a subquotient of 
$\GL_n(\Ff_\ell)$.
\item[b)] If $E/K$ is finite and $[E:K]$ is not divisible by $\ell$, then
$G_{\ell, E}^+=G_{\ell, K}^+$ and $EK_\ell^+=E_\ell^+$.
\end{enumerate}
\end{lemm} 

{\em Proof.} The profinite group $G_{\ell, E}$ is a closed normal subgroup
of $G_{\ell, K}$ and $G_{\ell, K}$ is isomorphic to a subquotient of $\GL_n(\Zz_\ell)$.
Hence there is a closed subgroup $U_\ell$ of $\GL_n(\Zz_\ell)$ and a closed normal
subgroup $V_\ell$ of $U_\ell$ such that there is an isomorphism $i: G_{\ell, E}\to U_\ell/V_\ell$. Furthermore there is a closed normal subgroup $U_\ell^+$ of $U_\ell$ containing
$V_\ell$ such that $i(G_{\ell, K}^+)=U_\ell^+/V_\ell$. The group $U_\ell/U_\ell^+$ is isomorphic to
$G_{\ell, E}/G_{\ell, E}^+$. Its order is coprime to $\ell$.  
The kernel of the restriction map 
$r\colon \GL_n(\Zz_\ell)\to \GL_n(\Ff_\ell)$ is a pro-$\ell$ group; hence the
intersection of this 
kernel with $U_\ell$ is contained in $U_\ell^+$. This shows that $r$ induces an
isomorphism  $U_\ell/U_\ell^+\to r(U_\ell)/r(U_\ell^+)$. Altogether we see that
$$\Gal(E_\ell^+/E)\cong G_{\ell, E}/G_{\ell, E}^+\cong U_\ell/U_\ell^+\cong r(U_\ell)/r(U_\ell^+) $$ 
and Part a) follows, 
because $r(U_\ell)/r(U_\ell^+) $ is obviously a subquotient of
$\GL_n(\Ff_\ell)$.

Every $\ell$-Sylow subgroup of $G_{\ell, E}$ lies in an $\ell$-Sylow subgroup
of $G_{\ell, K}$, hence $G_{\ell, E}^+\subset G_{\ell, K}^+$. Assume from
now on that $[E:K]$ is finite and not divisible by $\ell$. Then 
every $\ell$-Sylow subgroup of $G_{\ell, K}$ must map to the trivial group
under the projection $G_{\ell, K}\to G_{\ell, K}/G_{\ell, E}$, because
the order of the quotient group is coprime to $\ell$. Hence every 
$\ell$-Sylow subgroup of $G_{\ell, K}$ lies in $G_{\ell, E}$. This shows
that $G_{\ell, K}^+=G_{\ell, E}^+$. The Galois group
$\Gal(E_\ell/EK_\ell^+)$ is $G_{\ell, K}^+\cap G_{\ell, E}$ and the
Galois group $\Gal(E_\ell/EK_\ell^+)$ is $G_{\ell, E}^+$. As $G_{\ell, K}^+=G_{\ell, E}^+$ it follows
that $\Gal(E_\ell/EK_\ell^+)=\Gal(E_\ell/E_\ell^+)$, hence
$EK_\ell^+=E_\ell^+$.\hfill $\Box$

Let $S$ be a normal $\Qq$-variety with function field $K$. 
We shall now study families of representations of the fundamental group
$\pi_1(S)$ (viewing $S$ as a scheme equipped with the generic 
geometric point $\Spec(\widetilde{K})\to K$). 
Recall that we may identify $\pi_1(S)$ with
$\Gal(K_{S, \mathrm{nr}}/K)$.% where $K_{\mathrm{nr}}$ is the maximal extension of 
%$K$ in $\widetilde{K}$ which is unramified along $S$.

\begin{thm} Let $S/\Qq$ be a normal variety with function field $K$.
Let $P_\mathrm{nr}\in S_{\mathrm{nr}}$ be a closed point.
For every $\ell\in\Ll$ let $\Gamma_\ell$ be a profinite group and
$\rho_\ell\colon \pi_1(S)\to \Gamma_\ell$ a 
continuous homomorphism. We make two assumptions.
\begin{enumerate}
\item[a)] Assume there is an integer $n\in \Nn$ such that for every $\ell\in\Ll$ the profinite group $\Gamma_\ell$
is isomorphic to a subquotient of $\GL_n(\Zz_\ell)$.
\item[b)]
Assume that there exists an open subgroup $D'$ of the decomposition group $D_{K_{S, \mathrm{nr}}/K}(P_{\mathrm{nr}})$ such that the
family $(\rho_\ell)_{\ell\in\Ll}$ is independent over $D'$. 
\end{enumerate}
Then
the family $(\rho_\ell)_{\ell\in\Ll}$ is almost independent.\label{crit}
\end{thm}

This theorem may seem surprising at the first glance, since $D_{K_{S, \mathrm{nr}}/K}(P_{\mathrm{nr}})$ is usually far from being 
open in $\pi_1(S)$. The proof of Theorem \ref{crit} occupies the rest of this section. From now
on all the assumptions of Theorem \ref{crit} are in force, until the proof is finished.
For every algebraic extension $E/K$ contained in $K_{S, \mathrm{nr}}$ we define 
$G_{\ell, E}=\rho_\ell(\Gal(K_{S, \mathrm{nr}}/E))$, $G_{\ell, E}^+$, $E_\ell$ and $E_\ell^+$ exactly
as before. Furthermore we shall write $P_E$ for the point in
$S^{(E)}$ below $P_{\mathrm{nr}}$.

We tacitly assume in the sequal that $\widetilde{\Qq}$ denotes the algebraic closure
of $K$ {\em inside} $\widetilde{K}$. Then already $K_{S, \mathrm{nr}}$ contains
$\widetilde{\Qq}$, because the constant field extensions of $K$ are unramified along $S$.
The structure morphism $S_{\mathrm{nr}}\to\Spec(\Qq)$ factors through $\Spec(\widetilde{\Qq})$,
because $S_{\mathrm{nr}}$ is normal. It follows in particular that $k(P_{\mathrm{nr}})=
\widetilde{\Qq}$. 

\begin{lemm} There is a finite Galois extension 
$E/K$ contained in $K_{S, \mathrm{nr}}$ and a finite subset 
$I\subset \Ll$ such that the following 
statements about $E$ and $I$ hold true:
\begin{enumerate} 
\item[a)] For all $\ell\in\Ll\smallsetminus I$ the 
extension $E_\ell^+/E$ is a constant field extension, that is:
$\kappa_{E_\ell^+}E=E_\ell^+$. 
\item[b)] The point $P_E$ is a $\kappa_E$-rational point of $S^{(E)}$.
\item[c)] The family $(\rho_\ell)_{\ell\in \Ll\smallsetminus I}$ is independent
over $D_{K_{S, \mathrm{nr}}/E}(P_{\mathrm{nr}})$.\label{lemm2}
\end{enumerate}
\end{lemm}

{\em Proof.} Let $L:=\prod_{\ell\in\Ll} K_\ell^+$ be the composite field
of all the $K_\ell^+$. By Lemma \ref{lemm1}, for each $\ell\in\Ll$, the group 
$\Gal(K_\ell^+/K)$ is isomorphic to a subquotient
of $\GL_n(\Ff_\ell)$, and $|\Gal(K_\ell^+/K)|$ is not divisible by $\ell$. 
By \cite[Th\'eor\`eme 3']{bible} (which is a generalization 
due to Serre of the classical theorem of Jordan) it follows that there 
is an integer $d$ (independent of $\ell$) such that for every
$\ell\in \Ll$ the group $\Gal(K_\ell^+/K)$ has an abelian normal subgroup
$A_\ell$ of index $[\Gal(K_\ell^+/K):A_\ell]\le d$. Let $K_\ell'$ be the 
fixed field of $A_\ell$ in $K_\ell^+$. Then $K':=\prod_\ell K_\ell'$ is
a $d$-flat extension of $K$ and $K'K_\ell^+/K'$ is abelian for every $\ell\in\Ll$.
It follows that $L/K$ is a $d$-Jordanian extension. Furthermore 
$L/K$ is contained in $K_{S, \mathrm{nr}}$. By Proposition \ref{jordan},  
$L$ is a {\em finite} extension of $\kappa_LK$. Note that $\kappa_L/\Qq$ may well be
an infinite extension. Hence there is an element $\omega\in L$ such that $L=\kappa_LK(\omega)$.
Let $E_1$ be the Galois closure of $K(\omega)/K$ in $L$. Then $E_1/K$ is a finite
Galois extension and $\kappa_L E_1=L$. Hence we have a diagram of fields
$$\begin{xy}
  \xymatrix{
\kappa_L K\ar@{-}[r]\ar@{-}[d] & L\ar@{-}[d]\\
K\ar@{-}[r] & E_1
}
\end{xy}$$
in which the vertical extensions are constant field extensions and in which
the horizontal extensions are finite. Furthermore $L$ contains $K_\ell^+$ 
for every $\ell\in\Ll$. 

Now consider the canonical isomorphism
$$r\colon D_{K_{S, \mathrm{nr}}/K}(P_{\mathrm{nr}})\cong \Gal(k(P_{\mathrm{nr}})/{k(P_K)}).$$
Let $\lambda_1$ be the fixed field of $r(D')$ in $k(P_{\mathrm{nr}})=\widetilde{\Qq}$.
Since $D'$ is open in $D_{K_{S, \mathrm{nr}}/K}(P_{\mathrm{nr}})$, the field $\lambda_1$ is a 
finite extension of $k(P_K)$, so $\lambda_1$ is a finite extension of $\Qq$. 
Choose a finite Galois extension $\lambda/\kappa_K$ containing
$\lambda_1$ and $k(P_{E_1})$, and define $E:=\lambda E_1$. Then $S^{(E)}=S^{(E_1)}
\times_{\kappa_{E_1}} \Spec(\lambda)$ and $\kappa_E=\lambda$. There is the following
diagram of number fields:
$$\begin{xy}
  \xymatrix{
\lambda_1\ar@{-}[r]\ar@{-}[d] & \lambda\ar@{-}[d]\ar@{=}[r] & \kappa_E\ar@{-}[r] &k(P_E)\\
k(P_K)\ar@{-}[r]\ar@{-}[d] & k(P_{E_1})\ar@{-}[d]&&\\
\kappa_K\ar@{-}[r]& \kappa_{E_1}&&\\
}
\end{xy}$$

The fibre of $P_{E_1}$
under the projection $S^{(E)}\to S^{(E_1)}$ is $\Spec(\kappa_E\otimes_{\kappa_{E_1}} k(P_{E_1}))$, and this fibre splits up into the coproduct of $[k(P_{E_1}):\kappa_{E_1}]$ 
many copies of $\Spec(\kappa_E)=\Spec(\lambda)$, because $\lambda/\kappa_E$ is Galois
and $\lambda\supset k(P_{E_1})$. Thus all points in $S^{(E)}$ over $P_{E_1}$ are
$\kappa_E$-rational. In particular $P_E$ is $\kappa_E$-rational.
 
It follows that $$r(D_{K_{S, \mathrm{nr}}/E}(P_{\mathrm{nr}}))=
\Gal(k(P_{\mathrm{nr}})/{k(P_E)})=\Gal(k(P_{\mathrm{nr}})/\kappa_E),$$ 
and this group is an open subgroup of 
$r(D')=\Gal(k(P_{\mathrm{nr}})/\lambda_1)$ because 
$\kappa_E$ is a finite extension of $\lambda_1$. Hence $D_{K_{S, \mathrm{nr}}/E}(P_{\mathrm{nr}})$
is an open subgroup of $D'$.

As $(\rho_\ell)_{\ell\in\Ll}$ is independent over $D'$ by one
of our assumptions, it follows from part a) of Remark
\ref{indeprema} that there is a finite subset $I'\subset \Ll$ such
that the family $(\rho_\ell)_{\ell\in\Ll\smallsetminus I'}$ is independent
over $D_{K_{S, \mathrm{nr}}/E}(P_{\mathrm{nr}})$. Finally $K_\ell^+E/E$ is a 
constant field extension,
because $K_\ell^+E$ is an intermediate field of 
$LE/E$ and $LE=\kappa_LE$ is a constant field extension of $E$
due to our construction. By Lemma \ref{lemm1} we see that 
$E_\ell^+=K_\ell^+E$ for all $\ell\in \Ll$ which do not divide
the index $[E:K]$. Hence assertions a), b) and c) follow, if we 
put $I:=I'\cup\{\ell\in\Ll: \ell\ \mbox{divides}\ [E:K]\}$.\hfill $\Box$

\begin{lemm} Let $E$ and $I$ be as in Lemma \ref{lemm2}. Let 
$s\ge 1$. Let $\ell_1<\cdots<\ell_{s+1}$ be some elements of 
$\Ll\smallsetminus I$. Then 
$E_{\ell_1}\cdots E_{\ell_s}\cap E_{\ell_{s+1}}$ is a regular
extension of $\kappa_E$ (i.e. the algebraic closure of $\Qq$ in 
$E_{\ell_1}\cdots E_{\ell_s}\cap E_{\ell_{s+1}}$ is $\kappa_E$).\label{lemm3}
\end{lemm}

{\em Proof.} The canonical isomorphism
$$r\colon D_{K_{S, \mathrm{nr}}/E}(P_{\mathrm{nr}})\cong \Gal(k(P_{\mathrm{nr}})/k(P_E))$$
induces by restriction an isomorphism $$D_{K_{S, \mathrm{nr}}/E_\ell}(P_{\mathrm{nr}})=
D_{K_{S, \mathrm{nr}}/E}(P_{\mathrm{nr}})\cap \Gal(K_{S, \mathrm{nr}}/E_\ell)\cong \Gal(k(P_{\mathrm{nr}})/k(P_{E_\ell}))$$
for every $\ell\in\Ll$. Hence $k(P_{E_\ell})$ is the fixed field in
$k(P_{\mathrm{nr}})$ of the kernel of $\rho_\ell\circ r^{-1}$. The
family $(\rho_\ell)_{\ell\in \Ll\smallsetminus I}$ is independent over
$D_{K_{S, \mathrm{nr}}/E}(P_{\mathrm{nr}})$ by Lemma \ref{lemm2}. Hence Remark \ref{indeprem} shows that 
$(k(P_{E_\ell}))_{\ell\in\Ll\smallsetminus I}$ is linearly disjoint over
$k(P_E)$. Define $F:=E_{\ell_1}\cdots E_{\ell_s}$. 
There is a diagram of
residue fields:
$$
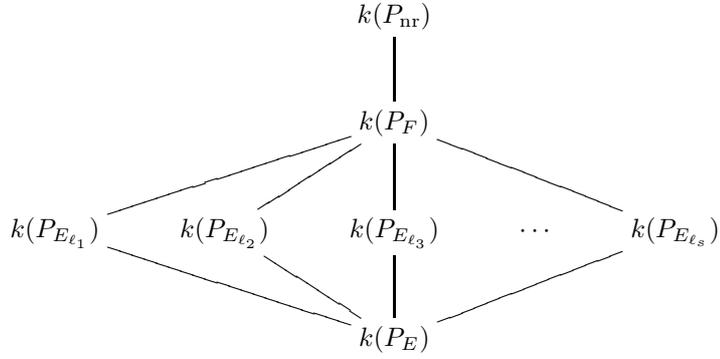
\begin{xy}
  \xymatrix{
 & & k(P_{\mathrm{nr}})\ar@{-}[d]\\ 
 & & k(P_F)\ar@{-}[lld]\ar@{-}[ld]\ar@{-}[d]\ar@{-}[rrd] & & \\
k(P_{E_{\ell_1}}) &k(P_{E_{\ell_2}}) & k(P_{E_{\ell_3}})  & \cdots & k(P_{E_{\ell_{s}}}) \\
 & & k(P_E)\ar@{-}[llu]\ar@{-}[lu]\ar@{-}[u]\ar@{-}[rru] & &
  }\end{xy}$$
We have $k(P_F)=k(P_{E_{\ell_1}})\cdots k(P_{E_{\ell_s}})$,
because 
$$\begin{array}{rcl}
\Gal(k(P_{\mathrm{nr}})/k(P_{E_{\ell_1}})\cdots k(P_{E_{\ell_s}}))&=&\bigcap_{i=1}^s G(k(P_{\mathrm{nr}})/k(P_{E_{\ell_i}}))=\\
&=& r(\bigcap_{i=1}^s D_{K_{S, \mathrm{nr}}/E_{\ell_i}}(P_{\mathrm{nr}}))=\\
&=&r(D_{K_{S, \mathrm{nr}}/E}(P_{\mathrm{nr}})\cap \bigcap_{i=1}^s \Gal(K_{S, \mathrm{nr}}/E_{\ell_i}))=\\
&=&r(D_{K_{S, \mathrm{nr}}/E}(P_{\mathrm{nr}})\cap \Gal(K_{S, {\mathrm{nr}}}/F))\\
&=& r(D_{K_{S, \mathrm{nr}}/F}(P_{\mathrm{nr}}))=G(k(P_{\mathrm{nr}})/k(P_F)).\end{array}$$ 
Furthermore there is a diagram
$$
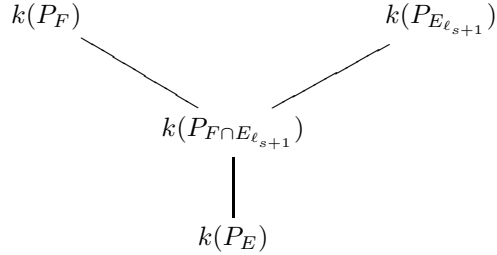
\begin{xy}
\xymatrix{
k(P_F)\ar@{-}[dr] & & k(P_{E_{\ell_{s+1}}})\ar@{-}[dl]\\
& k(P_{F\cap E_{\ell_{s+1}}})\ar@{-}[d] &\\
& k(P_E) &
}
\end{xy}$$ 
and $k(P_F)\cap k(P_{E_{\ell_{s+1}}})=k(P_E)$ due to the fact that $(k(P_{E_\ell}))_{\ell\in\Ll\smallsetminus I}$ is
linearly disjoint over $k(P_E)$. 
It follows that $k(P_{F\cap E_{\ell_{s+1}}})=k(P_E)$. Finally $k(P_E)=\kappa_E$, because $P_E$ is a $\kappa_E$-rational point
of $S^{(E)}$. This shows that the normalization of $S^{(E)}$ in 
$F\cap E_{\ell_{s+1}}$ has a $\kappa_E$-rational point and thus its function field $F\cap E_{\ell_{s+1}}$ must
be regular over $\kappa_E$. \hfill $\Box$

Let $\ell\ge 5$ be a prime number. We denote by $\Sigma_\ell$ the set of isomorphism classes 
of groups which are either the cyclic group $\Zz/\ell$, or the quotient of $\underline{H}(F)$ modulo
its center, where $F$ is a finite field of characteristic $\ell$ and $\underline{H}$ is a connected smooth
algebraic group over $F$ which is geometrically simple and simply connected. These are the simple groups of Lie
type in characteristic $\ell$. It is known (cf. \cite[Th\'eor\`eme 5]{bible}), that
$\Sigma_\ell\cap\Sigma_{\ell'}=\emptyset$ for all primes $5\le \ell<\ell'$. 
(As Serre points out in \cite{bible}, the proof of this theorem is essentially due to 
E. Artin \cite{artin}. It was completed in \cite{KLMS}.) In the following proof we shall strongly use
this result.

{\em End of Proof of Theorem \ref{crit}.} Let $E$ and $I$ be as in Lemma \ref{lemm2}.
In order to finish up the proof of Theorem
\ref{crit} it suffices to prove the following 

{\bf Claim.} {\em There is a finite subset $I'\subset \Ll$ containing $I$, such that 
$(E_\ell)_{\ell\in\Ll\smallsetminus I'}$ is linearly disjoint over $E$.}

In fact, once this claim is proven, it follows that $(\rho_\ell)_{\ell\in \Ll\smallsetminus I'}$ is
independent over $\Gal(K_{S, \mathrm{nr}}/E)$ by Remark \ref{indeprem}, and Remark \ref{indeprema} implies
that the whole family $(\rho_\ell)_{\ell\in \Ll}$ must be almost independent, as desired.

In \cite[Th\'eor\`eme 4]{bible} Serre proves: {\em There is
a constant $C$ such that for every prime number $\ell>C$ every
finite simple subquotient of $GL_n(\Zz_\ell)$ of order divisible by $\ell$ lies in $\Sigma_\ell$.}
This is a generalization of a well-known result of Nori (cf. \cite[Theorem B]{nori1987}).

Let us define $I':=I\cup \{2, 3\}\cup \{\ell\in \Ll: \ell\le C\}$. For $\ell\in\Ll$ every
non-trivial quotient of $G_{E, \ell}^+$ has order divisible by $\ell$: In fact, if
$h: G_{E, \ell}^+\to Q$ is an epimorphism onto a non-trivial group $Q$, 
then the image of some $\ell$-Sylow subgroup of $G_{E, \ell}$ under $h$ must be non-trivial.
Hence, for every $\ell\in \Ll\smallsetminus I'$ every finite
simple quotient of $G_{E, \ell}^+$ lies in $\Sigma_\ell$. 

We shall now prove the Claim. Let $s\ge 1$ and $\ell_1<\cdots<\ell_{s+1}$ be elements
of $\Ll\smallsetminus I'$. It suffices to show that $E_{\ell_1}\cdots E_{\ell_s}\cap E_{\ell_{s+1}}=E$,
assuming by induction that the sequence $(E_{\ell_1},\cdots, E_{\ell_s})$ is already linearly disjoint over $E$. 
This assumption implies 
$$\Gal(E_{\ell_1}\cdots E_{\ell_s}/E)\cong G_{\ell_1, E}\times \cdots \times G_{\ell_s, E}$$
and $\Gal(E_{\ell_1}\cdots E_{\ell_s}/E_{\ell_1}^+\cdots E_{\ell_s}^+)\cong G_{\ell_1, E}^+\times \cdots\times G_{\ell_s, E}^+$.
Suppose that $E_{\ell_1}\cdots E_{\ell_s}\cap E_{\ell_{s+1}}\neq E$. Then there would be an intermediate field $L$
of that extension such that $Q:=\Gal(L/E)$ is a finite simple group. We would have the following 
diagram of fields:
$$\begin{xy}
  \xymatrix{
  E_{\ell_1}\cdots E_{\ell_s}\ar@{-}[d]_{G_{\ell_1, E}^+\times \cdots \times G_{\ell_s, E}^+}\ar@{-}[dr] &    & E_{\ell_{s+1}}\ar@{-}[dl]\ar@{-}[d]^{G_{\ell_{s+1}, E}^+}\\
  E_{\ell_1}^+\cdots E_{\ell_s}^+\ar@{-}[dr]       & L\ar@{-}[d]^Q & E_{\ell_{s+1}}^+\ar@{-}[dl]\\
  & E &
}\end{xy}$$
But $L/\kappa_E$ is a regular extension (cf. Lemma \ref{lemm3}), hence $\kappa_L=\kappa_E$.
On the other hand $E_{\ell_i}^+/E$ is a constant field extension for every $i=1,\cdots, s+1$ (cf. Lemma \ref{lemm2}). 
It follows that 
$\Gal(LE_{\ell_{s+1}}^+/E_{\ell_{s+1}}^+)\cong Q$ and 
$\Gal(LE_{\ell_1}^+\cdots E_{\ell_s}^+/E_{\ell_1}^+\cdots E_{\ell_s}^+)\cong Q$. Hence $Q$ is simultaneously a quotient group
of $G_{\ell_1, E}^+\times\cdots\times G_{\ell_s, E}^+$ and of $G_{\ell_{s+1, E}}^+$. 
It follows that $$Q\in \left(\Sigma_{\ell_1}\cup\cdots\cup\Sigma_{\ell_s}\right)\cap\Sigma_{\ell_{s+1}},$$
which contradicts Artin's theorem that $\Sigma_\ell\cap\Sigma_{\ell'}=\emptyset$ for all primes 
$5\le \ell<\ell'$.\hfill $\Box$

\section{Proof of the main theorem}

{\em Proof of Theorem \ref{main2}.} Let $K$ be a finitely generated extension of $\Qq$. Let $X/K$ be
a separated scheme of finite type. 
Let $T=(T_1,\cdots, T_r)$ be a transcendence base of $K/\Qq$ and $S_0$ be 
the normalization of $\Spec(\Qq[T])$ in $K$. Then $S_0$ is a normal $\Qq$-variety with
function field $K$. 
The spreading-out principles in \cite{EGAIV3} (cf. in particular \cite{EGAIV3}[8.8.2],
\cite{EGAIV3}[8.10.5],
\cite{EGAIV3}[8.9.4]), allow
us to construct a dense open subscheme $S\subset S_0$ and a flat separated morphism of finite type
$f: \XX\to S$ with generic fibre $X$.  

We choose a closed point $P\in S$ and a closed point 
$P_{\mathrm{nr}}\in S_{\mathrm{nr}}$ over $P$ and denote by 
$\ol{P}: \Spec(k(P_{\mathrm{nr}}))\to S_{\mathrm{nr}}\to S$ the corresponding 
geometric point of $S$. Note that $k(P_{\mathrm{nr}})$ is algebraically closed (cf. the second paragraph after
Theorem \ref{crit}). We define $\widetilde{k}:=k(P_{\mathrm{nr}})$.
Furthermore we denote by $\ol{\xi}\colon \Spec(\widetilde{K})\to S$ 
the generic geometric point of $S$ afforded
by the choice of $\widetilde{K}$. We let $X_P:=\XX\times_S k(P)$, $X_{\ol{P}}=\XX\times_S \Spec(k(P_{\mathrm{nr}}))$ and 
$X_{\ol{\xi}}=\XX\times_S \Spec(\widetilde{K})$ be the corresponding fibres of 
$\XX$. Note that $X_{\ol{\xi}}=X_{\widetilde{K}}$ and $X_{\ol{P}}=X_{P, \widetilde{k}}$.

Let $q\in\Nn$. From now on we shall consider two cases. For the first case we define
$\rho_\ell:=\rho_{\ell, X}^{(q)}$, $T_\ell:=\Hrm^q(\Xx, \Zz_\ell)$, $V_\ell:=\Hrm^q(\Xx, \Qq_\ell)$,
$T_{\ell, P}:=\Hrm^q(\XP, \Zz_\ell)$, $V_{\ell, P}:=\Hrm^q(\XP, \Qq_\ell)$ and $\mathfrak{F}_\ell:=\Rrm^q f_*(\Zz_\ell)$ for every $\ell\in\Ll$. For the second case we define
$\rho_\ell:=\rho_{\ell, X, \mathrm{c}}^{(q)}$, $T_\ell:=\Hrm^q_\mathrm{c}(\Xx, \Zz_\ell)$, $V_\ell:=\Hrm^q_\mathrm{c}(\Xx, \Qq_\ell)$,
$T_{\ell, P}:=\Hrm^q_\mathrm{c}(\XP, \Zz_\ell)$, $V_{\ell, P}:=\Hrm^q_\mathrm{c}(\XP, \Qq_\ell)$ and $\mathfrak{F}_\ell:=\Rrm^q f_!(\Zz_\ell)$ for every $\ell\in\Ll$. 
In both cases $\rho_{\ell, P}$ will stand for the representation of $\Gal(\widetilde{k}/k(P))$ on $V_{\ell, P}$.

All residue characteristics of $S$ are zero. Hence  
there is a dense open subscheme $U\subset S$ such that for every $\ell\in\Ll$ the 
$\Zz_\ell$-sheaves $\Rrm^q f_*(\Zz_\ell)|U$ and  $\Rrm^q f_!(\Zz_\ell)|U$ are
lisse and of formation compatible with any base change $U'\to U$ (cf. \cite[Corollaire 2.6]{illusie}, \cite[Th\'eore\`me 3.1.2]{katzlaumon} and \cite[Th\'eore\`me 3.3.2]{katzlaumon}). Considering the cartesian diagrams
$$\begin{xy}
  \xymatrix{
  X_{\ol{P}}\ar[r]\ar[d] & \tilde{k}\ar[d] & \Xx\ar[r]\ar[d] & \widetilde{K}\ar[d]  \\
  f^{-1}(U)\ar[r] & U & f^{-1}(U)\ar[r] & U 
}\end{xy}$$
we can for every $\ell\in \Ll$ identify the stalks of $\mathfrak{F}_\ell$ by the following
base change isomorphisms
$$\mathfrak{F}_{\ell, \ol{P}}\cong T_{\ell, P}\ \ \mbox{and}\ \ 
\mathfrak{F}_{\ell, \ol{\xi}}\cong T_{\ell}.$$ 
The fact that the $\Zz_\ell$-sheaves $\mathfrak{F}_\ell|U$ are lisse  
implies that for every $\ell\in\Ll$ the representation $\rho_{\ell}$ 
factors through $\pi_1(U)$ and that 
there is a cospecialization isomorphism $\mathfrak{F}_{\ell, \ol{\xi}}\cong
\mathfrak{F}_{\ell, \ol{P}}$. Putting these isomorphisms together and tensoring with
$\Qq_\ell$ we obtain a cospecialization isomorphism 
$sp_{\ell}: V_\ell\cong V_{\ell, P}$
for every $\ell\in\Ll$. In order to take the Tate twists into account let
$\varepsilon_\ell: \Gal(\widetilde{K}/K)\to \Qq_\ell^\times$ be 
the cyclotomic character of $\Gal_K$ and by $\varepsilon_{\ell, P}: \Gal(\widetilde{k}/k(P))\to \Qq_\ell^\times$
 the cyclotomic character of $\Gal(\widetilde{k}/k(P))$. Let $d\in \Zz$ and define $\rho_\ell(d):=\rho_\ell\otimes
 \varepsilon_{\ell}^{\otimes d}$ and $\rho_{\ell, P}(d):=\rho_\ell\otimes
 \varepsilon_{\ell, P}^{\otimes d}$.
The cospecialization isomorphism $sp_\ell$ fits into a commutative
diagram
$$\begin{xy}
  \xymatrix{
 \Gal(K_{S, \mathrm{nr}}/K)\ar[r]^{\rho_\ell(d)} & \Aut_{\Qq_\ell}(V_\ell)\ar[dd] \\
 D_{K_{S, \mathrm{nr}}/K}(P_{\mathrm{nr}})\ar[u]\ar[d]  & \\
 \Gal(\widetilde{k}/k(P))\ar[r]^{\rho_{\ell, P}(d)} & \Aut_{\Qq_\ell}(V_{\ell, P})
}\end{xy}$$
for every $\ell\in\Ll$. 

There is a constant $b\in \Nn$ such that for every $\ell\in\Ll$
the inequality $\dim(V_\ell)\le b$ holds true (cf. \cite[Corollaire 1.3]{illusie}). 
Furthermore, if we donote the 
torsion part of the finitely generated $\Zz_\ell$-module
$T_\ell$ by $T_\ell'$, then $T_\ell/T_\ell'$ injects into $V_\ell$ and the representation $\rho_\ell(d)$
factors through $\Aut_{\Zz_\ell}(T_\ell/T_\ell')$. Hence $\mathrm{im}(\rho_\ell(d))$ (and also 
$\mathrm{im}(\rho_{\ell, P}(d))$) is isomorphic to
a closed subgroup of $\mathrm{GL}_b(\Zz_\ell)$ for every $\ell\in\Ll$. Hence the families 
$(\rho_\ell(d))_{\ell\in\Ll}$ and  
$(\rho_{\ell, P}(d))_{\ell\in\Ll}$ of representations of $\pi_1(U)$ satisfy assumption a) of 
Theorem \ref{crit} (and condition (B) of \cite[p.~3]{bible}). 

Now note that $X_P$ is a separated scheme of finite type over the number field $k:=k(P)$. For a place 
$v$ of a number field we denote by $p_v$ its residue characteristic. There is a finite 
extension $k'/k$ and a finite set $T$ of places of $k'$ such that the following holds true:
\begin{enumerate}
\item[(1)] For every place $v$ of $k'$ with $v\notin T$ and every $\ell\in \Ll\smallsetminus\{p_v\}$ the representation $\rho_{\ell, P}(d)$ is
unramified at $v$.
\item[(2)] For every $v\in T$, every extension $\hat{v}$ of $v$ to $\widetilde{k}$ and every $\ell\in \Ll\smallsetminus\{p_v\}$ the image of the inertia group $I_{\hat{v}}$ under the 
representation $\rho_{\ell, P}(d)$ is a pro-$\ell$ group.
\end{enumerate}
This is shown for $d=0$ in \cite[Th\'eore\`me 4.3]{illusie}, and the case $d\neq 0$ follows
as well, because the cyclotomic character $\varepsilon_{\ell, P}$ is unramified at every place $v$ of $k'$ 
with $p_v\neq \ell$. Because the family $(\rho_{\ell, P}(d))_{\ell\in\Ll}$ satisfies
the condition (B) of \cite[p.~3]{bible} and conditions (1) and (2) and because $k$ is a number field, 
Serre's theorem \cite[Th\'eore\`me 1]{bible} implies that the 
family $(\rho_{\ell, P}(d))_{\ell\in\Ll}$ is almost independent. 
Now the above diagram shows that there is an open subgroup $D'$ of $D_{K_{S, \mathrm{nr}}}(P_{\mathrm{un}})$ such that the restricted family $(\rho_{\ell}(d)|D')_{\ell\in\Ll}$ is independent, 
and our Theorem \ref{crit} implies that $(\rho_{\ell}(d))_{\ell\in\Ll}$ is almost independent as desired.\hfill $\Box$ 
\label{cohomology}

\appendix
\section{Abelian varieties}
The aim of this appendix is to give a more elementary direct proof of Corollary \ref{main1}, 
based on our independence 
criterion (cf. Theorem \ref{crit}) and on the corresponding 
results of Serre in the number field case. It avoids the use 
of \'etale cohomology.

{\em Proof of Corollary \ref{main1}.} Let $K$ be a finitely generated field of 
characteristic zero. Let $A/K$ be an abelian variety. It is enough to show that the 
family $(\eta_{\ell, A})_{\ell\in\Ll}$ defined in the introduction is almost independent. Then
Remark \ref{indeprem} implies the assertion.
There is
a normal $\Qq$-variety $S$ with function field $K$ and an  
abelian scheme $f: \mathcal{A}\to S$ with generic fibre $A$.

Let $P_{\mathrm{{{nr}}}}$ be a closed point of $S_{{\mathrm{nr}}}$ and $P$ the point of
$S$ below $P_{\mathrm{{{nr}}}}$. Then the residue field $k(P_{\mathrm{{{nr}}}})$ is an algebraic
closure of the number field $k(P)$. We define $\tilde{k}:=k(P_{\mathrm{{{nr}}}})$.
Let $A_P:=\mathcal{A}\times_S\Spec(k(P))$ be the special fibre of $\mathcal{A}$ at $P$. Then
$A_P$ is an abelian variety over the number field $k(P)$.

Let $n$ be an integer. The group scheme $\mathcal{A}[n]$ is finite and \'etale 
over $S$, because all residue characteristics of $S$ are zero.
Hence there is a finite extension $E/K$ contained in $K_{S, {\mathrm{nr}}}$ such that
$\mathcal{A}[n]\times_S S^{(E)}$ is a constant group scheme over $S^{(E)}$. In fact one can
take $E=K(A[n])$. This 
implies that both evaluation maps
$$\mathcal{A}[n](S^{(E)})\to A[n](E)\ \mbox{and}\ \mathcal{A}[n](S^{(E)})\to A_P[n](\tilde{k})$$
are isomorphisms. In particular the action of $\Gal_K$ on $A[n](\tilde{K})$ factors
through $\Gal(K_{S, {\mathrm{nr}}}/K)$ (and in fact through $\Gal(E/K)$). We obtain a composite
isomorphism
$$A[n](\tilde{K})\cong \mathcal{A}[n](S^{(E)})\cong A_P[n](\tilde{k}).$$
Taking limits, we obtain for each $\ell\in\Ll$ an isomorphism
$$T_\ell(A)\cong T_\ell(A_P)$$
and the action of $\Gal_K$ on $T_\ell(A)$ factors through $\Gal(K_{S, {\mathrm{nr}}}/K)$.
This isomorphism fits into a commutative diagram
$$\begin{xy}
  \xymatrix{
 \Gal(K_{S, {\mathrm{nr}}}/K)\ar[r]^{\eta_{\ell, A}} & \Aut_{\Zz_\ell}(T_\ell(A))\ar[dd]\\
 D_{K_{S, \mathrm{nr}}/K}(P_{\mathrm{nr}})\ar[u]\ar[d]  & \\
 \Gal(k(P_{\mathrm{{{nr}}}})/k(P))\ar[r]^{\eta_{\ell, A_P}} & \Aut_{\Zz_\ell}(T_\ell(A_P))
}\end{xy}$$
Recall that $A_P$ is an abelian variety over the number field $k(P)$. 
Hence Serre's theorem (cf. \cite[Section 3]{bible}) implies that 
the family $(\eta_{\ell, A_P})_{\ell\in\Ll}$ is almost independent. It follows
that there is an open subgroup $D'$ in $D_{K_{S, \mathrm{nr}}/K}(P_{\mathrm{nr}})$ such that 
the family $(\eta_{\ell, A})_{\ell\in\Ll}$ is independent over $D'$. 
Now, by our Theorem \ref{crit}, the family 
$(\eta_{\ell, A})_{\ell\in\Ll}$ must be almost independent, as desired. \hfill $\Box$

\begin{center}
{\bf Acknowledgements}
\end{center}

We want to thank Luc Illusie for a detailed list of comments on an earlier version of this paper. We 
are indebted to him for pointing out to us how to apply results from \cite{illusie} and \cite{katzlaumon}
in order to remove
an unnecessary smoothness assumption from our main Theorem.  
Furthermore we want to thank Cornelius Greither and Moshe Jarden for encouragement and discussions, and for
a variety of very helpful comments.  
Both authors were supported by the Deutsche Forschungsgemeinschaft research
grant GR 998/5-1. Sebastian Petersen thanks Adam Mickiewicz University in Pozna\'n for its hospitality 
during several research visits and Wojciech Gajda thanks Universit\"at der Bundeswehr in Munich for its
hospitality during a visit in January 2011. In Pozna\'n Wojciech Gajda was partially supported by a reasearch
grant of the Polish Ministry of Science and Higher Education. The mathematical content of the present work has been much 
influenced by the preprint \cite{bible} of Serre, and also by the inspiring article \cite{katzlang} of
Katz and Lang. We acknowledge this with pleasure.

%\appendix
%\section{Abelian varieties}
%The aim of this appendix is to give a direct proof of Corollary \ref{main1} which is independent of \'etale 
%cohomology.

{\sc Wojciech Gajda\\
Faculty of Mathematics and Computer Science\\
Adam Mickiewicz University\\
Umultowska 87\\
61614 Pozna\'{n}, Poland}\\
E-mail adress: \texttt{\small gajda@amu.edu.pl}
\par\medskip

{\sc Sebastian Petersen\\
FB 10 - Mathematik und Naturwissenschaften\\ 
Universit\"at Kassel\\
Heinrich-Plett-Str. 40\\
34132 Kassel, Germany}\\
E-mail adress: \texttt{\small petersen@mathematik.uni-kassel.de} 
\end{document}